\documentclass[a4paper]{amsart}

\usepackage{graphicx}
\usepackage{amsthm,amsmath,amsfonts,amssymb,mathrsfs}
\usepackage{lmodern}  
\usepackage{bbm} 
\usepackage{enumitem}

\newtheorem{definition}{Definition}[section]

\newtheorem{proposition}{Proposition}[section]

\newtheorem{remark}{Remark}[section]
\newtheorem{notation}{Notation}[section]

\newcommand{\df}[1]{\textbf{\textit{#1}}}  
\newcommand{\Tau}{\mathcal{T}} 
\newcommand{\E}{\mathbbmss{E}}
\newcommand{\U}{\mathbbmss{U}}
\newcommand{\PP}{\mathbb{P}} 
\newcommand{\Acal}{\mathcal{A}}  
\newcommand{\Bcal}{\mathcal{B}}  
\newcommand{\Scal}{\mathcal{S}}  
\newcommand{\SSE}[1][\U]{{\mathcal{SS}(#1)}_\E}  
\newcommand{\SSG}[2]{{\mathcal{SS}(#1)}_{#2}} 
\newcommand{\SPE}[1][\U]{{\mathcal{SP}(#1)}_\E}  
\newcommand{\softsubseteq}{\tilde{\subseteq}}  
\newcommand{\softequal}{\tilde{=}}  
\newcommand{\softin}{\tilde{\in}}  
\newcommand{\softnotequal}{\tilde{\ne}}  
\newcommand{\softnotin}{\tilde{\notin}}  
\newcommand{\nullsoftset}{(\tilde{\emptyset},\E)}  
\newcommand{\absolutesoftset}[1][\U]{(\tilde{#1},\E)}  
\newcommand{\absolutesoftsetG}[2]{(\tilde{#1},#2)}  
\newcommand{\softsetminus}{\widetilde{\setminus}}  
\newcommand{\softcap}{\tilde{\cap}}  
\newcommand{\softbigcup}{\widetilde{\bigcup}}  
\newcommand{\softbigcap}{\widetilde{\bigcap}}  
\newcommand{\softcl}[2][X]{s\textit{-}cl_{#1}\left( #2 \right)}  
\newcommand{\softhomeomorphic}{\tilde{\thickapprox}}  

\begin{document}
\title{An Embedding Lemma in Soft Topological Spaces}
\author{Giorgio Nordo}
\date{}
\maketitle

\begin{abstract}
In 1999, Molodtsov initiated the concept of Soft Sets Theory as a new mathematical tool and a completely different
approach for dealing with uncertainties in many fields of applied sciences.
In 2011, Shabir and Naz introduced and studied the theory of soft topological spaces, also defining
and investigating many new soft properties as generalization of the classical ones.
In this paper, we introduce the notions of soft separation between soft points and soft closed sets in order
to obtain a generalization of the well-known Embedding Lemma for soft topological spaces.
\end{abstract}

\section{Introduction}
Many practical problems in economics, engineering, environment, social science,
medical science etc. cannot be studied by classical methods, because they
have inherent difficulties due to the inadequacy of the theories of parameterization tools
in dealing with uncertainties.
In 1999, Molodtsov \cite{molodtsov} initiated the novel concept
of Soft Sets Theory as a new mathematical tool
and a completely different approach for dealing with uncertainties
while modelling problems in computer science, engineering physics, economics,
social sciences and medical sciences.
Molodtsov defines a soft set as a parameterized family of subsets of universe set where each element
is considered as a set of approximate elements of the soft set.

In 2011, Shabir and Naz \cite{shabir} introduced the concept of soft topological spaces,
also defining and investigating the notions of soft closed sets, soft closure,
soft neighborhood, soft subspace and some separation axioms.
Some other properties related to soft topology were studied by
\c{C}a\u{g}man, Karata\c{s} and Enginoglu in \cite{cagman2011}.
In the same year Hussain and Ahmad \cite{hussain} continued the study investigating
the properties of soft closed sets, soft neighbourhoods, sof interior, soft exterior
and soft boundary.

In the present paper we will present the notions of family of soft mappings
soft separating soft points and soft points from soft closed sets
in order to give a generalization of the well-known Embedding Lemma for soft topological spaces.

\section{Preliminaries}
In this section we present some basic definitions and results of the theories
of soft sets and soft topological spaces, simplifying them in a suitable way whenever possible.
Terms and undefined concepts are used as in \cite{engelking}. 


\begin{definition}{\rm\cite{molodtsov}}
\label{def:softset}
Let $\U$ be an initial universe set and $\E$ be a nonempty set of parameters (or abstract attributes)
under consideration with respect to $\U$ and $A\subseteq \E$,
we say that a pair $(F,A)$ is a \df{soft set} over $\U$
if $F$ is a set-valued mapping $F: A \to \PP(\U)$
which maps every parameter $e \in A$ to a subset $F(e)$ of $\U$.
\end{definition}

In other words, a soft set is not a real (crisp) set
but a parameterized family $\left\{ F(e) \right\}_{e\in A}$ of subsets of the universe $\U$.
For every parameter $e \in A$, $F(e)$ may be considered as the set of \textit{$e$-approximate elements}
of the soft set $(F,A)$.

\begin{remark}
\label{rem:sameparameter}
In 2010, Ma, Yang and Hu \cite{ma} proved that every soft set $(F,A)$ is equivalent
to the soft set $(F,\E)$ related to the whole set of parameters $\E$,
simply considering empty every approximations of parameters which are missing in $A$,
that is extending in a trivial way its set-valued mapping,
i.e. setting $F(e)=\emptyset$, for every $e \in \E \setminus A$.
\\
For such a reason, in this paper we can consider all the soft sets over the same parameter set $\E$
as in \cite{chiney} and we will redefine all the basic operations and relations
between soft sets originally introduced in \cite{molodtsov,maji2002,maji2003} as in \cite{nazmul},
that is by considering the same parameter set.
\end{remark}

\begin{definition}{\rm\cite{zorlutuna}}
\label{def:setofsoftsets}
The set of all the soft sets over a universe $\U$ with respect to a set of parameters $\E$
will be denoted by $\SSE$.
\end{definition}

\begin{definition}{\rm\cite{nazmul}}
\label{def:softsubset}
Let $(F,\E),(G,\E) \in \SSE$ be two soft sets over a common universe $\U$
and a common set of parameters $\E$,
we say that $(F,\E)$ is a \df{soft subset} of $(G,\E)$ and we write
$(F,\E) \softsubseteq (G,\E)$
if $F(e)\subseteq G(e)$ for every $e \in \E$.
\end{definition}

\begin{definition}{\rm\cite{nazmul}}
\label{def:softequal}
Let $(F,\E),(G,\E) \in \SSE$ be two soft sets over a common universe $\U$, we say that
$(F,\E)$ and $(G,\E)$ are \df{soft equal} and we write $(F,\E) \softequal (G,\E)$
if $(F,\E) \softsubseteq (G,\E)$ and $(G,\E) \softsubseteq (F,\E)$.
\end{definition}

\begin{definition}{\rm\cite{nazmul}}
\label{def:nullsoftset}
A soft set $(F,\E)$ over a universe $\U$ is said to be \df{null soft set}
and it is denoted by $\nullsoftset$ if $F(e) = \emptyset$ for every $e \in \E$.
\end{definition}

\begin{definition}{\rm\cite{nazmul}}
\label{def:absolutesoftset}
A soft set $(F,\E) \in \SSE$ over a universe $\U$ is said to be a \df{absolute soft set}
and it is denoted by $\absolutesoftset$
if $F(e) = \U$ for every $e \in \E$.
\end{definition}

\begin{definition}{\rm\cite{nazmul}}
\label{def:softcomplement}
Let $(F,\E) \in \SSE$ be a soft set over a universe $\U$, the \df{soft complement}
(or more exactly the \textit{soft relative complement}) of $(F,\E)$,
denoted by $(F,\E)^\complement$, is the soft set $\left( F^\complement, E \right)$
where $F^\complement : \E \to \PP(\U)$ is the set-valued mapping
defined by $F^\complement(e) = F(e)^\complement = \U \setminus F(e)$, for every $e \in \E$.
\end{definition}

\begin{definition}{\rm\cite{nazmul}}
\label{def:softdifference}
Let $(F,\E),(G,\E) \in \SSE$ be two soft sets over a common universe $\U$,
the \df{soft difference} of $(F,\E)$ and $(G,\E)$,
denoted by $(F,\E) \softsetminus (G,\E)$, is the soft set $\left( F \setminus G, E \right)$
where $F \setminus G : \E \to \PP(\U)$ is the set-valued mapping
defined by $(F \setminus G)(e) = F(e) \setminus G(e)$, for every $e \in \E$.
\end{definition}

Clearly, for every soft set $(F,\E) \in \SSE$, it results
$(F,\E)^\complement \, \softequal \, \absolutesoftset \softsetminus (F,\E)$.

\begin{definition}{\rm\cite{nazmul}}
\label{def:generalizedsoftunion}
Let $\left\{(F_i,\E) \right\}_{i\in I} \subseteq \SSE$ be a nonempty subfamily
of soft sets over a universe $\U$,
the (generalized) \df{soft union} of $\left\{(F_i,\E) \right\}_{i\in I}$,
denoted by $\softbigcup_{i \in I} (F_i,\E) $,
is defined by $\left(\bigcup_{i \in I} F_i, \E \right)$
where $\bigcup_{i \in I} F_i: \E \to \PP(\U)$ is the set-valued mapping
defined by $\left(\bigcup_{i \in I} F_i\right)(e) = \bigcup_{i \in I} F_i(e)$, for every $e \in \E$.
\end{definition}

\begin{definition}{\rm\cite{nazmul}}
\label{def:generalizedsoftintersection}
Let $\left\{(F_i,\E) \right\}_{i\in I} \subseteq \SSE$ be a nonempty subfamily
of soft sets over a universe $\U$,
the (generalized) \df{soft intersection} of $\left\{(F_i,\E) \right\}_{i\in I}$,
denoted by $\softbigcap_{i \in I} (F_i,\E) $,
is defined by $\left(\bigcap_{i \in I} F_i, E \right)$
where $\bigcap_{i \in I} F_i: \E \to \PP(\U)$ is the set-valued mapping
defined by $\left(\bigcap_{i \in I} F_i\right)(e) = \bigcap_{i \in I} F_i(e)$, for every $e \in \E$.
\end{definition}

\begin{definition}{\rm\cite{al-khafaj}}
\label{def:softdisjunct}
Two soft sets $(F,\E)$ and $(G,\E)$ over a common universe $\U$
are said to be \df{soft disjoint} if their soft intersection is the soft null set,
i.e. if $(F,\E) \softcap (G,\E) \, \softequal \, \nullsoftset$.
\end{definition}

\begin{definition}{\rm\cite{xie}}
\label{def:softpoint}
A soft set $(F,\E) \in \SSE$ over a universe $\U$ is said to be a \df{soft point} over $U$
if it has only one non-empty approximation and it is a singleton,
i.e. if there exists some parameter $\alpha \in E$
and an element $p \in \U$ such that
$F(\alpha) = \{ p \}$ and $F(e)=\emptyset$ for every $e \in E \setminus \{ \alpha \}$.
Such a soft point is usually denoted by$(p_\alpha, \E)$.
The singleton $\{ p \}$ is called the \textit{support set} of the soft point
and $\alpha$ is called the \textit{expressive parameter} of $(p_\alpha, \E)$.
\end{definition}

\begin{definition}{\rm\cite{xie}}
\label{def:setofsoftpoints}
The set of all the soft points over a universe $\U$ with respect to a set of parameters $\E$
will be denoted by $\SPE$.
\end{definition}

Since any soft point is a particular soft set, it is evident that $\SPE \subseteq \SSE$.

\begin{definition}{\rm\cite{xie}}
\label{def:softpointsoftbelongstosoftset}
Let $(p_\alpha, \E) \in \SPE$ and $(F,\E) \in \SSE$
respectively be a soft point and a softset over a common universe $\U$.
We say that \df{the soft point $(p_\alpha, \E)$ soft belongs to the soft set $(F,\E)$}
and we write $(p_\alpha, \E) \softin (F,\E)$, if the soft point is a soft subset of the soft set,
i.e. if $(p_\alpha, \E) \softsubseteq (F,\E)$
and hence if $p \in F(\alpha)$.
\end{definition}

\begin{definition}{\rm\cite{das}}
\label{def:equalitysoftpoints}
Let $(p_\alpha, \E), (q_\beta, \E) \in \SPE$ be two soft points over a common universe $\U$,
we say that $(p_\alpha, \E)$ and $(q_\beta, \E)$ are \df{soft equal},
and we write $(p_\alpha, \E) \softequal (q_\beta, \E)$,
if they are equals as soft sets and hence if $p = q$ and $\alpha = \beta$.
\end{definition}

\begin{definition}{\rm\cite{das}}
\label{def:distinctssoftpoints}
We say that two soft points $(p_\alpha, \E)$ and $(q_\beta, \E)$ are \df{soft distincts},
and we write $(p_\alpha, \E) \softnotequal (q_\beta, \E)$,
if and only if $p\ne q$ or $\alpha \ne \beta$.
\end{definition}

According to Remark \ref{rem:sameparameter} the following definitions by
Kharal and Ahmad have been simplified and slightly modified for soft sets on a common parameter set.

\begin{definition}{\rm\cite{kharal}}
\label{def:softmapping}
Let $\SSG{\U}{\E}$ and $\SSG{\U'}{\E'}$ be two sets of soft open sets
over the universe sets $\U$ and $\U'$ with respect to the sets of parameters $\E$ and $\E'$, respectively.
and consider a mapping $\varphi: \U \to \U'$ between the two universe sets and
a mapping $\psi: \E \to \E'$ between the two set of parameters.
The mapping $\varphi_\psi : \SSG{\U}{\E} \to \SSG{\U'}{\E'}$
which maps every soft set $(F,\E)$ of $\SSG{\U}{\E}$
to a soft set $\varphi_\psi \left( (F,\E) \right)$ of $\SSG{\U'}{\E'}$ denoted by $(\varphi_\psi(F), \E')$
where $\varphi_\psi(F) : \E' \to \PP(\U')$ is the set-valued mapping
defined by $\varphi_\psi(F)(e') = \bigcup_{e\in \psi^{-1}(\{e'\})} \varphi(F(e))$
for every $e' \in \E'$,
is called a \df{soft mapping} from $\U$ to $\U'$ induced by the mappings $\varphi$ and $\psi$,
while the soft set $\varphi_\psi (F,\E) \softequal (\varphi_\psi(F), \E')$
is said to be the \df{soft image} of the soft set $(F,\E)$
under the soft mapping $\varphi_\psi : \SSG{\U}{\E} \to \SSG{\U'}{\E'}$.
\\
The soft mapping $\varphi_\psi : \SSG{\U}{\E} \to \SSG{\U'}{\E'}$ is said
\df{injective} (respectively \df{surjective}, \df{bijective})
if the mappings $\varphi: \U \to \U'$ and $\psi: \E \to \E'$ are both
injective (resp. surjective, bijective).
\end{definition}

\begin{definition}{\rm\cite{kharal}}
\label{def:softinverseimage}
Let $\varphi_\psi : \SSG{\U}{\E} \to \SSG{\U'}{\E'}$ be a soft mapping
induced by the mappings $\varphi: \U \to \U'$ and $\psi: \E \to \E'$
between the two sets $\SSG{\U}{\E}, \SSG{\U'}{\E'}$ of soft sets
and consider a soft set $(G,\E')$ of $\SSG{\U'}{\E'}$.
The \df{soft inverse image} of $(G,\E')$ under the soft mapping
$\varphi_\psi : \SSG{\U}{\E} \to \SSG{\U'}{\E'}$,
denoted by $\varphi_\psi^{-1} \left( (G,\E') \right)$ is the
soft set $(\varphi_\psi^{-1}(G), \E')$ of $\SSG{\U}{\E}$
where $\varphi_\psi^{-1}(G) : \E \to \PP(\U)$ is the set-valued mapping
defined by $\varphi_\psi^{-1}(G)(e) = \varphi^{-1}\left( G\left( \psi(e) \right) \right)$
for every $e \in \E$.
\end{definition}

The concept of soft topological space as topological space defined by a family of soft sets
over a initial universe with a fixed set of parameters
was introduced in 2011 by Shabir and Naz \cite{shabir}.

\begin{definition}{\rm\cite{shabir}}
\label{def:softtopology}
Let $X$ be an initial universe set, $\E$ be a nonempty set of parameters with respect to $X$
and $\Tau \subseteq \SSE[X]$ be a family of soft sets over $X$, we say that
$\Tau$ is a \df{soft topology} on $X$ with respect to $\E$ if the following four conditions are satisfied:
\begin{enumerate}[label=(\roman*)]
\item the null soft set belongs to $\Tau$, i.e. $\nullsoftset \in \Tau$
\item the absolute soft set belongs to $\Tau$, i.e. $\absolutesoftset[X] \in \Tau$
\item the soft intersection of any two soft sets of $\Tau$ belongs to $\Tau$, i.e.
for every $(F,\E), (G,\E) \in \Tau$ then $(F,\E) \softcap (G,\E) \in \Tau$.
\item the soft union of any subfamily of soft sets in $\Tau$ belongs to $\Tau$, i.e.
for every $\left\{(F_i,\E) \right\}_{i\in I} \subseteq \Tau$ then
$\softbigcup_{i \in I} (F_i,\E) \in \Tau$
\end{enumerate}
The triplet $(X, \Tau, \E)$ is called a \df{soft topological space} over $X$ with respect to $\E$.
\\
In some case, when it is necessary to better specify the universal set and the set of parameters,
the topology will be denoted by $\Tau(X,\E)$.
\end{definition}

\begin{definition}{\rm\cite{shabir}}
\label{def:softopenset}
Let $(X,\Tau,\E)$ be a soft topological space over $X$ with respect to $\E$,
then the members of $\Tau$ are said to be \df{soft open set} in $X$.
\end{definition}

\begin{definition}{\rm\cite{hazra}}
\label{def:comparisonofsofttopologies}
Let $\Tau_1$ and $\Tau_2$ be two soft topologies over a common universe set $X$
with respect to a set of paramters $\E$.
We say that $\Tau_2$ is \df{finer} (or stronger) than $\Tau_1$
if $\Tau_1 \subseteq \Tau_2$ where $\subseteq$ is the usual set-theoretic relation of inclusion between crisp sets.
In the same situation, we also say that $\Tau_1$ is \df{coarser} (or weaker) than $\Tau_2
$.
\end{definition}

\begin{definition}{\rm\cite{shabir}}
\label{def:softclosedset}
Let $(X,\Tau,\E)$ be a soft topological space over $X$ and be $(F,\E)$ be a soft set over $X$.
We say that $(F,\E)$ is \df{soft closed set} in $X$ if its complement $(F,\E)^\complement$
is a soft open set, i.e. if $(F^\complement,\E) \in \Tau$.
\end{definition}

\begin{notation}
The family of all soft closed sets of a soft topological space $(X,\Tau,\E)$ over $X$ with respect to $\E$
will be denoted by $\sigma$,
or more precisely with $\sigma(X,\E)$ when it is necessary to specify
the universal set $X$ and the set of parameters $\E$.
\end{notation}

\begin{definition}{\rm\cite{aygunoglu}}
\label{def:softbase}
Let $(X,\Tau,\E)$ be a soft topological space over $X$ and $\Bcal \subseteq \Tau$ be a non-empty subset of soft open sets.
We say that $\Bcal$ is a \df{soft open base} for $(X,\Tau,\E)$ if every soft open set of $\Tau$ can be
expressed as soft union of a subfamily of $\Bcal$, i.e. if for every $(F,\E) \in \Tau$
there exists some $\Acal \subset \Bcal$ such that $(F,\E) = \softbigcap \left\{ (A,\E) : \, (A,\E) \in \Acal \right\}$.
\end{definition}

\begin{definition}{\rm\cite{zorlutuna}}
\label{pro:softneighbourhoodofasoftpoint}
Let $(X, \Tau, \E)$ be a soft topological space, $(N,\E) \in \SSE[X]$ be a soft set
and $(x_\alpha, \E) \in \SPE[X]$ be a soft point over a common universe $X$.
We say that $(N,\E)$ is a \df{soft neighbourhood} of the soft point $(x_\alpha, \E)$
if there is some soft open set soft containing the soft point and soft contained in the soft set,
that is if there exists some soft open set $(A,\E) \in \Tau$
such that $(x_\alpha, \E) \softin (A,\E) \softsubseteq (N,\E)$.
\end{definition}

\begin{definition}{\rm\cite{shabir}}
\label{def:softclosure}
Let $(X,\Tau,\E)$ be a soft topological space over $X$ and $(F,\E)$ be a soft set over $X$.
Then the \df{soft closure} of the soft set $(F,\E)$,
denoted by $\softcl{(F,\E)}$,
is the soft intersection of all soft closed set over $X$ soft containing $(F,\E)$, that is
$$\softcl{(F,\E)} \softequal \, \softbigcap \left\{ (C,\E) \in \sigma(X,\E) :
\, (F,\E) \softsubseteq (C,\E) \right\}$$
\end{definition}

\begin{definition}{\rm\cite{zorlutuna}}
\label{def:softcontinuousmapping}
Let $\varphi_\psi : \SSG{X}{\E} \to \SSG{X'}{\E'}$ be a soft mapping
between two soft topological spaces $(X,\Tau,\E)$ and $(X',\Tau',\E')$
induced by the mappings $\varphi: X \to X'$ and $\psi: \E \to \E'$
and $(x_\alpha, \E) \in \SPE[X]$ be a soft point over $X$.
We say that the soft mapping $\varphi_\psi$ is \df{soft continuous at the soft point $(x_\alpha, \E)$}
if for each soft neighbourhood $(G,\E')$
of $\varphi_\psi \left( (x_\alpha, \E) \right)$ in $(X',\Tau',\E')$
there exists some soft neighbourhood $(F,\E)$ of $(x_\alpha, \E)$ in $(X,\Tau,\E)$
such that $\varphi_\psi \left( (F,\E) \right) \softsubseteq (G,\E')$.
\\
If $\varphi_\psi$ is soft continuous at every soft point $(x_\alpha, \E) \in \SPE[X]$,
then $\varphi_\psi : \SSG{X}{\E} \to \SSG{X'}{\E'}$  is called \df{soft continuous} on $X$.
\end{definition}

\begin{definition}{\rm\cite{aygunoglu}}
\label{def:softsubbase}
Let $(X,\Tau,\E)$ be a soft topological space over $X$ and $\Scal \subseteq \Tau$ be
a non-empty subset of soft open sets.
We say that $\Scal$ is a \df{soft open subbase} for $(X,\Tau,\E)$ if the family of all finite soft intersections
of members of $\Scal$ forms a soft open base for $(X,\Tau,\E)$.
\end{definition}

\begin{proposition}{\rm\cite{aygunoglu}}
\label{pro:softtopologyinducedbyasoftsubbase}
Let $\Scal \subseteq \SSG{X}{\E}$ be a a family of soft sets over $X$ and containing both
the null soft set $\nullsoftset$ and the absolute soft set $\absolutesoftset[X]$.
Then the family $\Tau(\Scal)$
of all soft union of finite soft intersections of soft sets in $\Scal$
is a soft topology having $\Scal$ as subbase.
\end{proposition}

\begin{definition}{\rm\cite{aygunoglu}}
\label{def:softtopologyinducedbyasoftsubbase}
Let $\Scal \subseteq \SSG{X}{\E}$ be a a family of soft sets over $X$ respect to a set of parameters $\E$
and such that $\nullsoftset, \absolutesoftset[X] \in \Scal$, then
the soft topology $\Tau(\Scal)$ of the above Proposition \ref{pro:softtopologyinducedbyasoftsubbase}
is called the \df{soft topology generated} by the soft subbase $\Scal$ over $X$
and $(X,\Tau(\Scal), \E)$ is said to be the \df{soft topological space generated by $\Scal$}.
\end{definition}

\section{An Embedding Lemma for soft topological spaces}

\begin{definition}{\rm\cite{aygunoglu}}
\label{def:softinitialtopology}
Let $\SSG{X}{\E}$ be the set of soft sets over a universe set $X$
with respect to a set of parameter $\E$
and consider a family of soft topological spaces $\left\{ (Y_i, \Tau_i, \E_i ) \right\}_{i \in I}$
and a corresponding family $\left\{ (\varphi_\psi)_i \right\}_{i\in I}$
of soft mappings $(\varphi_\psi)_i : \SSG{X}{\E} \to \SSG{Y_i}{\E_i}$.
Then the soft topology $\Tau(\Scal)$ generated by the soft subbase
$\Scal = \left\{ (\varphi_\psi)_i^{-1} \left( (G,\E_i) \right) : \,
(G,\E_i) \softin \Tau_i , \, i \in I \right\}$
of all soft inverse images of the soft mappings $(\varphi_\psi)_i$
is called the \df{initial soft topology}
induced on $X$ by the family of soft mappings $\left\{ (\varphi_\psi)_i \right\}_{i\in I}$
and it is denoted by $\Tau_{ini} \left(X, \E, Y_i, \E_i, (\varphi_\psi)_i; i\in I \right)$ .
\end{definition}

\begin{definition}
\label{def:softprojectionmappings}
Let $\left\{ (X_i, \Tau_i, \E_i ) \right\}_{i \in I}$ be a family of soft topological spaces
over the universe sets $X_i$ with respect to the sets of parameters $\E_i$, respectively.
For every $i \in I$, the soft mapping
${\left(\pi_i \right)}_{\rho_i} : \SSG{\prod_{i \in I} X_i}{\prod_{i \in I} \E_i} \to \SSG{X_i}{\E_i}$
induced by the classical projections $\pi_i : \prod_{i \in I} X_i \to X_i$
and $\rho_i : \prod_{i \in I} \E_i \to \E_i$
is said the \df{$i$-th soft projection mapping}
and, by setting $(\pi_\rho)_i = {\left(\pi_i \right)}_{\rho_i}$, it will be denoted by
$(\pi_\rho)_i : \SSG{\prod_{i \in I} X_i}{\prod_{i \in I} \E_i} \to \SSG{X_i}{\E_i}$.
\end{definition}

\begin{definition}{\rm\cite{aygunoglu}}
\label{def:softtopologicalproduct}
Let $\left\{ (X_i, \Tau_i, \E_i ) \right\}_{i \in I}$ be a family of soft topological spaces
and consider the corresponding family $\left\{ (\pi_\rho)_i \right\}_{i \in I}$
of soft projection mappings
$(\pi_\rho)_i : \SSG{\prod_{i \in I} X_i}{\prod_{i \in I} \E_i} \to \SSG{X_i}{\E_i}$
(with $i \in I$).
Then, the initial soft topology
$\Tau_{ini} \left( \prod_{i \in I} X_i, \E, X_i, \E_i, (\pi_\rho)_i; i\in I \right)$
induced on $\prod_{i \in I} X_i$ by the family of soft projection mappings
$\left\{ (\pi_\rho)_i \right\}_{i \in I}$
is called the \df{soft topological product} of the soft topological space
$( X_i, \Tau_i, \E_i )$ (with $i \in I$)
and denoted by $\Tau \left( \prod_{i \in I} X_i \right)$.
\end{definition}

\begin{definition}{\rm\cite{aras}}
\label{def:softhomeomorphism}
Let $(X,\Tau,\E)$ and $(X',\Tau',\E')$ be two soft topological spaces
over the universe sets $X$ and $X'$ with respect to the sets of parameters $\E$ and $\E'$, respectively.
We say that a soft mapping $\varphi_\psi : \SSG{X}{\E} \to \SSG{X'}{\E'}$
is a \df{soft homeomorphism} if it is soft continuous, bijective
and its inverse $\varphi_\psi^{-1} : \SSG{X'}{\E'} : \to \SSG{X}{\E} $
is a soft continuous mapping too.
In such a case, the soft topological spaces $(X,\Tau,\E)$ and $(X',\Tau',\E')$
are said \df{soft homeomorphic} and we write that $(X,\Tau,\E) \softhomeomorphic (X',\Tau',\E')$.
\end{definition}

\begin{definition}
\label{def:softcembedding}
Let $(X,\Tau,\E)$ and $(X',\Tau',\E')$ be two soft topological spaces.
We say that a soft mapping $\varphi_\psi : \SSG{X}{\E} \to \SSG{X'}{\E'}$
is a \df{soft embedding} if its corestriction
$\varphi_\psi : \SSG{X}{\E} \to \varphi_\psi \left( \SSG{X}{\E} \right)$
is a soft homeomorphism.
\end{definition}

\begin{definition}
\label{def:softdiagonalmapping}
Let $(X,\Tau,\E)$ be a soft topological space over a universe set $X$ with respect to a set of parameter $\E$,
let $\left\{ (X_i, \Tau_i, \E_i ) \right\}_{i \in I}$ be a family of soft topological spaces
over a universe set $X_i$ with respect to a set of parameters $\E_i$, respectively
and consider a family $\left\{ (\varphi_\psi)_i \right\}_{i\in I}$
of soft mappings $(\varphi_\psi)_i = {\left( \varphi_i \right)}_{\psi_i}: \SSG{X}{\E} \to \SSG{X_i}{\E_i}$
induced by the mappings $\varphi_i : X \to X_i$ and $\psi_i : \E \to \E_i$ (with $i \in I$).
Then the soft mapping
$\Delta = \varphi_\psi : \SSG{X}{\E} \to \SSG{\prod_{i \in I} X_i}{\prod_{i \in I} \E_i} $
induced by the diagonal mappings (in the classical meaning)
$\varphi = \Delta_{i \in I} \varphi_i : X \to \prod_{i \in I} X_i$
and $\psi  = \Delta_{i \in I} \psi_i: \E \to \prod_{i \in I} \E_i$
(respectively defined by $\varphi(x) = \langle \varphi_i(x) \rangle_{i \in I}$ for every $x \in X$
and by $\psi(e) = \langle \psi_i(e) \rangle_{i \in I}$ for every $x \in X$)
is called the \df{soft diagonal mapping} of the soft mappings $(\varphi_\psi)_i$ (with $i\in I$)
and denoted by $\Delta = \Delta_{i \in I} (\varphi_\psi)_i : \SSG{X}{\E} \to \SSG{\prod_{i \in I} X_i}{\prod_{i \in I} \E_i} $.
\end{definition}

\begin{definition}
\label{def:familyofsofttmappingseparatingsoftpoints}
Let $\left\{ (\varphi_\psi)_i \right\}_{i\in I}$ be a family of
of soft mappings $(\varphi_\psi)_i : \SSG{X}{\E} \to \SSG{X_i}{\E_i}$
between a soft topological space $(X,\Tau,\E)$
and a family of soft topological spaces $\left\{ (X_i, \Tau_i, \E_i ) \right\}_{i \in I}$.
We say that the family $\left\{ (\varphi_\psi)_i \right\}_{i\in I}$
\df{soft separates soft points} of $(X,\Tau,\E)$
if for every $(x_\alpha, \E), (y_\beta, \E) ( \in \SPE[X]$ such that
$(x_\alpha, \E) \softnotequal (y_\alpha, \E)$ there exists some $i \in I$ such that
$(\varphi_\psi)_i (x_\alpha, \E) \softnotequal (\varphi_\psi)_i (y_\beta, \E)$.
\end{definition}

\begin{definition}
\label{def:familyofsofttmappingseparatingsoftpointsfromsoftclosedsets}
Let $\left\{ (\varphi_\psi)_i \right\}_{i\in I}$ be a family of
of soft mappings $(\varphi_\psi)_i : \SSG{X}{\E} \to \SSG{X_i}{\E_i}$
between a soft topological space $(X,\Tau,\E)$
and a family of soft topological spaces $\left\{ (X_i, \Tau_i, \E_i ) \right\}_{i \in I}$.
We say that the family $\left\{ (\varphi_\psi)_i \right\}_{i\in I}$
\df{soft separates soft points from soft closed sets} of $(X,\Tau,\E)$
if for every $(C,\E) \in \sigma(X,\E)$
and every $(x_\alpha, \E) ( \in \SPE[X]$ such that
$(x_\alpha, \E) \softin \absolutesoftsetG{X}{\E} \softsetminus (C,\E)$
there exists some $i \in I$ such that
$(\varphi_\psi)_i (x_\alpha, \E) \softnotin \softcl[X_i]{(\varphi_\psi)_i (C,\E) }$.
\end{definition}

\begin{proposition}[\textbf{Soft Embedding Lemma}]
\label{pro:softembeddinglemma}
Let $(X,\Tau,\E)$ be a soft topological space,
$\left\{ (X_i, \Tau_i, \E_i ) \right\}_{i \in I}$ be a family of soft topological spaces
and $\left\{ (\varphi_\psi)_i \right\}_{i\in I}$ be a family of
of soft continuous mappings $(\varphi_\psi)_i : \SSG{X}{\E} \to \SSG{X_i}{\E_i}$
that separates both the soft points
and the soft points from the soft closed sets of $(X,\Tau,\E)$.
Then the diagonal mapping $\Delta = \Delta_{i \in I} (\varphi_\psi)_i : \SSG{X}{\E}
\to \SSG{\prod_{i \in I} X_i}{\prod_{i \in I} \E_i} $
of the soft mappings $(\varphi_\psi)_i$ is a soft embedding.
\end{proposition}

\section{Conclusion}
In this short announcement paper we have introduced the notions of family of soft mappings
separating points and points from closed sets 
and that of soft diagonal mapping
in order to define the necessary framework
for proving a generalization to soft topological spaces of the well-known Embedding Lemma
for classical (crisp) topological spaces.
Such a result could be the start point for investigating other important topics
in soft topology such as extension and compactifications theorems,
metrization theorems etc.
Details and proofs will be given in a next extended paper.


\vskip 6mm
\parskip=0pt {\noindent {\sc Giorgio NORDO \newline
MIFT - Dipartimento di Scienze Matematiche e Informatiche, scienze Fisiche e scienze della Terra, \newline
Viale F. Stagno D'Alcontres, 31 --
Contrada Papardo, salita Sperone, 31 - 98166 Sant'Agata -- Messina
(ITALY)}} \vskip 1pt \noindent
E-mail:  {\tt giorgio.nordo@unime.it}


\begin{thebibliography}{99}

\bibitem{aras} Aras C.G., Sonmez A., \c{C}akalli H. 2013, On soft mappings,
{\it Proceedings of CMMSE 2013 - $13^{th}$ International Conference on Computational and Mathematical Methods in Science and Engineering},
arXiv:1305.4545, 11 pages. 

\bibitem{aygunoglu} Ayg\"{u}no\u{g}lu A., Ayg\"{u}n H. 2012, Some notes on soft topology spaces,
{\it Neural Computing and Applications} {\bf 21}, pp. 113-119.

\bibitem{cagman2011} \c{C}a\u{g}man N., Karata\c{s} S., Enginoglu S. 2011, Soft topology,
{\it Computers \& Mathematics with Applications} {\bf 62}, pp. 351-358.

\bibitem{chiney} Chiney M., Samanta S.K. 2017, Soft topology redefined,
arXiv:1701.00466, 18 pages. 

\bibitem{das} Das S., Samanta S.K. 2013, Soft metric,
{\it Annals of Fuzzy Mathematics and Informatics} {\bf 6} (1), pp. 77-94.

\bibitem{engelking} Engelking R. 1989, {\it General Topology}
(Berlin: Heldermann Verlag).

\bibitem{hazra} Hazra H., Majumdar P., Samanta S.K. 2012, Soft topology,
{\it Fuzzy information and Engineering} {\bf 4} (1), pp. 105-115.

\bibitem{hussain} Hussain S., Ahmad B. 2011, Some properties of soft topological spaces,
{\it Computers \& Mathematics with Applications} {\bf 62} (11), pp. 4058-4067.

\bibitem{al-khafaj} Al-Khafaj M.A.K., Mahmood M.H. 2014,
Some Properties of Soft Connected Spaces and Soft Locally Connected Spaces,
{\it IOSR Journal of Mathematics} {\bf 10} (5), pp. 102-107.

\bibitem{kharal} Kharal A., Ahmad B., 2011, Mappings on soft classes,
{\it New Mathematics and Natural Computation} {\bf 7} (3), pp. 471-481.

\bibitem{ma} Ma Z., Yang B., Hu B. 2010, Soft set theory based on its extension,
{\it Fuzzy Information and Engineering} {\bf 2} (4), pp. 423-432.

\bibitem{maji2002} Maji, P.K., Roy A.R. 2002, An application of soft sets in a decision making problem,
{\it Computers \& Mathematics with Applications} {\bf 44}, pp. 1077-1083.

\bibitem{maji2003} Maji P.K., Biswas R., Roy A.R. 2003, Soft set theory,
{\it Computers \& Mathematics with Applications} {\bf 45}, pp. 555-562.

\bibitem{molodtsov} Molodtsov D. 1999, Soft set theory -- first result,
{\it Computers \& Mathematics with Applications} {\bf 37}, pp. 19-31.

\bibitem{nazmul} Nazmul S., Samanta S.K. 2013, Neighbourhood properties of soft topological spaces,
{\it Annals of Fuzzy Mathematics and Informatics} {\bf 6} (1), pp. 1-15.

\bibitem{shabir} Shabir M., Naz M. 2011, On soft topological spaces,
{\it Computers \& Mathematics with Applications} {\bf 61}, pp. 1786-1799.

\bibitem{xie} Xie N. 2015, Soft points and the structure of soft topological spaces,
{\it Annals of Fuzzy Mathematics and Informatics} {\bf 10} (2), pp. 309-322.

\bibitem{zorlutuna} Zorlutuna \.{I}, Akdag M., Min W.K., Atmaca S. 2012, Remarks on soft topological spaces,
{\it Annals of Fuzzy Mathematics and Informatics} {\bf 3} (2), pp. 171-185.

\end{thebibliography}
\end{document}